\documentclass[10pt]{article}
\usepackage{latexsym}\usepackage{amsbsy}\usepackage{amssymb}
\usepackage[latin1]{inputenc} 

\newcommand{\ep}{\hspace*{\fill}$\Box$}
\newcommand{\eps}{\varepsilon}
\newcommand{\pr}{\noindent{\bf Proof. }}
\newcommand{\ms}{\medskip\\}

\newcommand{\R}{\mathbb R}
\newcommand{\N}{\mathbb N}
\newcommand{\C}{\mathbb C}

\newcommand{\K}{\mathbb K}

\newcommand{\gs}{\ensuremath{{\mathcal G}} }

\newcommand{\es}{\ensuremath{{\mathcal E}} }
\newcommand{\esm}{\ensuremath{{\mathcal E}_M} }

\newcommand{\ns}{\ensuremath{{\mathcal N}} }
\newcommand{\ks}{\ensuremath{{\mathcal K}} }

\newcommand{\comp}{\subset\subset}
\newcommand{\cinfty}{{\cal C}^\infty}
\newtheorem{thr}{\hspace*{-3mm} \bf}[section]
\newcommand{\bt}{\begin{thr} {\bf Theorem. }}
\newcommand{\et}{\end{thr}}
\newcommand{\bp}{\begin{thr} {\bf Proposition. }}
\newcommand{\bc}{\begin{thr} {\bf Corollary. }}
\newcommand{\blem}{\begin{thr} {\bf Lemma. }}
\newcommand{\bex}{\begin{thr} {\bf Example. }\rm} 
\newcommand{\bexs}{\begin{thr} {\bf Examples. }\rm}
\newcommand{\bd}{\begin{thr} {\bf Definition. }}
\newcommand{\beast}{\begin{eqnarray*}}
\newcommand{\eeast}{\end{eqnarray*}}

\newcommand{\al}{\alpha}
\newcommand{\bet}{\beta} 

\newcommand{\Om}{\Omega}\newcommand{\Ga}{\Gamma}

\newcommand{\vphi}{\varphi}

\newcommand{\D}{{\cal D}}\newcommand{\Vol}{\mbox{Vol\,}}

\newcommand{\CC}{{\cal C}}

\newcommand{\nn}{\nonumber}
\newcommand{\beq}{ \begin{equation} }\newcommand{\eeq}{\end{equation} }
\newcommand{\bea}{\begin{eqnarray}}\newcommand{\eea}{\end{eqnarray}}
\newcommand{\beas}{\begin{eqnarray*}}\newcommand{\eeas}{\end{eqnarray*}}
\newcommand{\beqs}{\begin{equation*}}\newcommand{\eeqs}{\end{equation*}}

\newcommand{\GL}{\mbox{GL}}\newcommand{\bfs}{\boldsymbol}
\newcommand{\ben}{\begin{enumerate}}\newcommand{\een}{\end{enumerate}}
\newcommand{\ba}{\begin{array}}\newcommand{\ea}{\end{array}}
\newcommand{\brem}{\begin{thr} {\bf Remark. }\rm}
\newcommand{\ethi}{\end{thr}}
\newcommand{\pro}{\mathrm{pr}}
\newcommand{\esmh}{\ensuremath{{\mathcal E}^{\mathit h}_M} }
\newcommand{\gsh}{\ensuremath{{\mathcal G}^{\mathit h}} }
\newcommand{\simh}{\ensuremath{\sim_{\mathit h}}}

\begin{document}

\title{Intrinsic characterization of manifold-valued generalized functions}
\author{Michael Kunzinger \footnote{Electronic mail: michael.kunzinger@univie.ac.at}\\ 
        Roland Steinbauer \footnote{Electronic mail: roland.steinbauer@univie.ac.at}\\
         {\small Department of Mathematics, University of Vienna}\\
         {\small Strudlhofg.\ 4, A-1090 Wien, Austria}\\
        James A.\ Vickers \footnote{Electronic mail: J.A.Vickers@maths.soton.ac.uk}\\
         {\small  University of Southampton, Faculty of Mathematical Studies,}\\ 
         {\small Highfield, Southampton SO17 1BJ, United Kingdom}
       }
\date{April 2002}  
\maketitle
                                                                           
\begin{abstract}
The concept of generalized functions taking values
in a differentiable manifold (\cite{gfvm, curvature}) is extended to 
a functorial theory. We establish several characterization results which
allow a global intrinsic formulation both of the theory of manifold-valued generalized functions
and of generalized vector bundle homomorphisms.
As a consequence, a characterization of equivalence that does not resort to derivatives 
(as provided in \cite{found} for the scalar-valued cases of Colombeau's construction) is achieved. 
These results are employed to derive a point value description of all types of generalized functions
valued in manifolds and to show that composition can be carried out unrestrictedly.
Finally, a new concept of association adapted to the present setting is introduced.
 
\end{abstract}

\section{Introduction}\label{intro}
While originating as a tool in the field of nonlinear PDEs
a growing number of applications of nonlinear generalized functions
\cite{c1,c2,MObook} in a geometrical context---especially in the theory of 
general relativity (cf., e.g., \cite{clarke,hb5,wilson}, as well as \cite{vickersESI}
for a survey) and Lie group analysis of differential equations 
(cf.\ \cite{symm,DKP,book})---has triggered the development of a {\em geometric} 
theory of algebras of generalized functions (\cite{book}). 
Based on 
\cite{AB, CM, Jel},
a diffeomorphism invariant scalar theory 
in the so called ``full''
setting (which allows for a distinguished embedding of the space
of Schwartz distributions) was developed in \cite{found,vim}.
On the other hand the task of ``geometrizing'' the special setting of
Colombeau's construction---which will also provide the framework for 
the present article---was 
started 
in \cite{RD,DP,ndg,curvature}. 
In particular, in \cite{ndg} a theory of generalized sections of vector bundles
showing maximal compatibility with the distributional setting of 
\cite{deR,marsden} was introduced. Also a point value description of
generalized functions on differentiable manifolds was achieved.
This setting in turn was used to set up a ``generalized (pseudo-)Riemannian 
geometry'' in \cite{curvature}. Finally, in \cite{gfvm}
the study of generalized functions {\rm taking values} in a differentiable
manifold was initiated. More precisely, the space
$\gs[X,Y]$ of generalized functions on the manifold $X$ and
taking values in the manifold $Y$ was defined as was the space
$\mathrm{Hom}_{\gs}(E,F)$ of generalized vector bundle homomorphisms
from $E$ to $F$. When dealing with generalized mappings
on manifolds the need for such concepts arises quite naturally;
for example when considering the flow of a generalized vector field as well as 
geodesics of a generalized metric.

The aim of the present paper is to complete this
approach by extending it to a functorial theory. In particular, we
provide several global characterizations of the notions of moderateness
and equivalence (which replaces the notion of negligibility in the absence of
a linear structure) for generalized functions from $X$ to $Y$ (Section \ref{emch}) 
as well as for generalized vector bundle homomorphisms from $E$ to $F$ (Section \ref{homch}). 
These characterizations in turn enable us to give a positive answer to 
the question raised in \cite{gfvm}, Remark 2.11 whether equivalence
is equivalent to equivalence of order zero. More precisely, we establish
also in the case of both $\gs[X,Y]$ and $\mathrm{Hom}_{\gs}(E,F)$
a characterization of equivalence of nets  without making any reference to 
derivatives provided the nets are assumed to be moderate (Theorems 
\ref{equivchar} and \ref{equivvbchar}). Analogous results have 
been proven for virtually all variants of scalar- or vector valued Colombeau generalized 
functions (cf.\ \cite{found}, Th.\ 13.1, and \cite{ndg}, Sec.\ 4).
In the context of $\gs[X,Y]$ and $\mathrm{Hom}_{\gs}(E,F)$ 
these characterizations provide the key to proving that
(a) generalized functions taking values in a manifold as well
as generalized vector bundle homomorphisms are characterized by their 
generalized point values (Theorem \ref{pvchar} respectively \ref{pvcharhom}), and
(b) composition of generalized functions between manifolds as well as
of generalized vector bundle homomorphisms can be carried out unrestrictedly
(Theorems \ref{generalcomp} and \ref{pvvbth}).
Analogous results are given in Section \ref{hybch} for ``hybrid'' 
Colombeau spaces $\gs^h[X,F]$ of generalized maps defined on a manifold $X$ 
and taking values in a vector bundle $F$ over $Y$. These spaces were 
introduced in \cite{curvature}, Sec.\ 4 to allow for the notion of generalized
sections along generalized maps and, in particular, to define the notion of
a geodesic of a generalized pseudo-Riemannian metric. Moreover we prove 
that for any $u\in\gs[X,Y]$ the space of generalized vector bundle homomorphisms 
on $u$, $\mathrm{Hom}_u(E,F)$ can be endowed with a natural vector space structure. 
While trivial in the smooth setting the corresponding 
result in the present framework
relies on the fact that for any $v\in\mathrm{Hom}_u(E,F)$ we may
construct a representative whose induced map on the base coincides with a given
representative of $u$ (Proposition 4.5).
Finally in Section \ref{coupledcal} we consider the concept of  
`coupled calculus' and, in particular, the notion 
of zero-association which in the context of $\gs[X,Y]$ 
allows us to make statements regarding compatibility of this approach 
with classical global analysis. We start this work by recalling
the basic notions from the geometric theory of generalized functions 
(in the special
setting) and, in particular, of generalized functions valued in a smooth 
manifold. 

\section{Geometric theory of generalized functions}\label{geomth} 

We begin by fixing some notation. Throughout this paper $X$ and $Y$ denote 
smooth paracompact Hausdorff manifolds of dimension $n$ and $m$ respectively. 
Vector bundles with base space $X$ will be denoted 
by $(E,X,\pi)$ (or $(E,X,\pi_X)$) and for a chart $(V,\vphi)$ in $X$, 
a vector bundle chart $(V,\Phi)$ over $\varphi$ will be written in the form 
($\K=\R$ respectively $\C$)
\begin{equation} \label{vbchart}
{\displaystyle
\begin{array}{rcl}
      \Phi:\,\pi^{-1}(V)&\to&\vphi(V)\times\K^{n'}\nn\\
        z&\to&(\vphi(p),\vphi^1(z),\dots,\vphi^{n'}(z))
        \equiv(\vphi(p),\bfs{\vphi}(z))\,, 
\end{array}
}
\end{equation}
where $p=\pi(z)$ and the typical fiber is $\K^{n'}$.
Given a vector bundle atlas $(V_\al,\Phi_\al)_\al$ we write 
$\Phi_\al\circ\Phi_\beta^{-1}(y,w) = (\vphi_{\al\bet}(y),\bfs{\vphi_{\al\beta}}
(y)w)$, where $\vphi_{\al\bet}:=\vphi_\al\circ\vphi^{-1}_\bet$ is the change 
of chart on the base and $\bfs{\vphi_{\al\beta}}:$ $\vphi_\bet(V_\al\cap V_\beta)
\to \GL({n'},\K)$ denotes the transition functions.
The space of smooth sections of the vector bundle $(E,X,\pi)$ 
is denoted by  $\Ga(X,E)$, the space of smooth ($r,s$)-tensor fields (i.e.,
$E=T^r_s$) by ${\cal T}^r_s$ and the space of smooth vector bundle homomorphisms 
from $E$ to $F$ by  $\mathrm{Hom}(E,F)$.
If $f\in \mathrm{Hom}(E$, $F)$ we write $\underline{f}: X\to Y$ for the smooth 
map induced on bases, i.e., $\pi_Y\circ f = \underline{f}\circ \pi_X$.
Local vector bundle homomorphisms with respect to 
vector bundle charts $(V,\Phi)$ of $E$ and $(W,\Psi)$ of $F$, i.e.,
$f_{\mathrm{\Psi}\mathrm{\Phi}}:=
\mathrm{\Psi}\circ f \circ \mathrm{\Phi}^{-1}: \vphi(V\cap \underline{f}^{-1}(W))
\times \K^{n'} \to \psi(W) \times \K^{m'}$
will be written in the form
\beas
f_{\mathrm{\Psi}\mathrm{\Phi}}(x,\xi) = 
(f_{\mathrm{\Psi}\mathrm{\Phi}}^{(1)}(x),f_{\mathrm{\Psi}\mathrm{\Phi}}^{(2)}(x)
\cdot\xi)\,.
\eeas

Turning now to notions from Colombeau theory we set $I=(0,1]$ and ${\mathcal E}(X)=\cinfty(X)^I$
and let ${\cal P}(X)$ denote the space of linear differential operators on $X$. We define the
spaces of moderate and  negligible nets in ${\mathcal E}(X)$ by
\beas
        \esm(X)&:=&\{ (u_\eps)_\eps\in{\mathcal E}(X):\ 
        \forall K\subset\subset X,\ \forall P\in{\cal P}(X)\ \exists N\in\N:
  \\&&\hphantom{(u_\eps)_\eps\in{\mathcal E}(X):\ \forall K\subset\subset X,\ \forall 
     P\in{\cal P}(\}} 
        \sup_{p\in K}|Pu_\eps(p)|=O(\eps^{-N})\}
   \\
        \ns(X)&:=& \{ (u_\eps)_\eps\in\esm(X):\ 
        \forall K\subset\subset X,\ \forall q \in\N_0:\
        \sup_{p\in K}|u_\eps(p)|=O(\eps^{q}))\}\,.
\eeas
$\gs(X):= \esm(X)/\ns(X)$ is called the special\footnote{Since we are going 
to work entirely in the ``special'' setting of Colombeau's construction we 
omit this term henceforth.} Colombeau algebra on $X$ and we denote its
elements by $u=[(u_\eps)_\eps)]$. $\gs(\_)$ is a fine sheaf of differential algebras
with respect to the Lie derivative along smooth vector fields defined by 
$L_\xi u=[(L_\xi u_\eps)_\eps]$. $u$ is in $\gs(X)$ if and only 
if $u\circ\psi_\al\in\gs(\psi_\al(V_\al))$ (the local Colombeau algebra on 
$\psi_\al(V_\al)$) for all charts $(V_\al,\psi_\al)$. 
$\cinfty(X)$ is a subalgebra of $\gs(X)$ and there exist injective sheaf morphisms 
embedding $\D'(X)$, the space of Schwartz distributions on $X$, into $\gs(X)$.

A net $(p_\eps)_\eps \in X^I$ of points in $X$ is called compactly supported 
if $p_\eps$ stays in a fixed compact set for small $\eps$; the set of generalized points
is denoted by
$X_c$. Two nets $(p_\eps),\, (q_\eps)_\eps\in X^I$ are called equivalent, $(p_\eps)_\eps\sim (q_\eps)_\eps$, 
if $d_h(p_\eps,q_\eps) = O(\eps^m)$ for each $m>0$, where $d_h$ denotes the distance 
function induced on $X$ by any Riemannian metric $h$. The quotient space $\widetilde X_c$ 
of the set of compactly supported points modulo $\sim$ is called the space of
compactly supported generalized points on $X$ and we write $\tilde p = [(p_\eps)_\eps]$. 
Inserting $\tilde p$ into  $u\in \gs(X)$ yields a well-defined element $[(u_\eps(p_\eps))_\eps]$
of $\ks$, the space of generalized numbers 
(corresponding to $\K=\R$ respectively $\C$ and defined
as the set of moderate nets of numbers $(r_\eps)_\eps \in \K^I$ with
$|r_\eps| = O(\eps^{-N})$ for some $N$ modulo negligible nets
$|r_\eps| = O(\eps^{m})$ for each $m$). Moreover, $u\in \gs(X)$ is uniquely 
determined by its point values on $\widetilde X_c$, i.e., $u=v\ \Leftrightarrow
u(\tilde p)=v(\tilde p)\ \forall \tilde p\in\widetilde X_c$ (\cite{ndg}, Th.\ 1).

The $\gs(X)$-module $\Ga_\gs(X,E)$ of generalized sections in $E$ is defined using 
analogous asymptotic estimates with respect to the norm on the fibers induced by any 
Riemannian metric on $X$. Setting $\Gamma_{\mathcal E}(X,E)=\Ga(X,E)^I$ 
and denoting by ${\mathcal P}(X,E)$ the space of linear differential operators on $\Ga(X,E)$
we define 
\beas
        \Gamma_{\esm}(X,E)&:=& \{ (s_\eps)_{\eps}\in \Gamma_{\mathcal E}(X,E) :
                \ \forall P\in {\mathcal P}(X,E)\, \forall K\comp X \, \exists N\in \N:\\
                 &&\hspace{4cm}\sup_{p\in K}\|Pu_\eps(p)\| = O(\eps^{-N})\}\\
        \Gamma_\ns(X,E)&:=& \{ (s_\eps)_{\eps}\in \Gamma_{\esm}(X,E) :
                \ \forall K\comp X \, \forall m\in \N:\\
                 &&\hspace{4.5cm}\sup_{p\in K}\|u_\eps(p)\| = O(\eps^{m})\}\,,\
\eeas
Finally $\Gamma_\gs(X,E):=\Gamma_{\esm}(X,E)/\Gamma_\ns(X,E)$.
Any $s=[(s_\eps)_\eps]\in\Gamma_\gs(X,E)$ corresponds to a family 
$(s_\al)_\al=((s^i_\al)_\al)_{i=1}^{n'}$,
where $s_\al$ is called the local expression of $s$ with its components
$s^i_\al:=\Psi^i_\al\circ s\circ\psi_\al^{-1} \in\gs(\psi_\al(V_\al))$
($i=1,\dots, {n'}$) satisfying $s^i_\al(x)\,=\,(\bfs{\psi}_{\al\beta})^i_j
(\psi_\beta\circ\psi^{-1}_\al(x))\,s^j_\beta$ $
(\psi_\beta\circ\psi^{-1}_\al(x))$ for all $x\in \psi_\al(V_\al\cap V_\beta)$.
$\Ga_\gs(\_\,,E)$ is a fine sheaf of  projective and finitely generated  
$\gs(X)$-modules (\cite{ndg}, Th.\ 5) and we have 
\begin{equation}\label{tensorp}
  \Gamma_\gs(X,E)=\gs(X)\otimes\Ga(X,E)\,,
\end{equation}
where the tensor product is taken over $\CC^\infty(X)$ (\cite{ndg}, Th.\ 4).

For later reference we now explicitly recall the main definitions from \cite{gfvm}.
The space $\gs[X,Y]$ of compactly bounded (c-bounded) 
generalized Co\-lom\-beau functions from $X$ to $Y$
is defined as the quotient of the set of $\es_M[X,Y]$ 
of moderate, c-bounded maps from $X$ to $Y$ by a certain equivalence 
relation defined below.
\bd \label{modmapmf}
$\esm[X,Y]$ is defined as the set of all $(u_\eps)_\eps \in \cinfty(X,Y)^I$ 
satisfying
\begin{itemize}
  \item[(i)] $\forall K\comp\Om\ \exists \eps_0>0\  \exists K'\comp Y \  \forall
             \eps<\eps_0:\ u_\eps(K) \subseteq K'$.
  \item[(ii)]
   $\forall k\in\N$, for each chart $(V,\vphi)$ in $X$, each 
   chart $(W,\psi)$ in $Y$, each $L\comp V$ and each $L'\comp W$
   there exists $N\in \N$  with
   $$\sup\limits_{x\in L\cap u_\eps^{-1}(L')} \|D^{(k)}
   (\psi\circ u_\eps \circ \vphi^{-1})(\vphi(p))\| =O(\eps^{-N}).
   $$
\end{itemize}
\et

\bd\label{equrel} $(u_\eps)_\eps$ and $(v_\eps)_\eps \in \esm[X,Y]$ are
called equivalent, $(u_\eps)_\eps \sim (v_\eps)_\eps$,
if the following conditions are satisfied:
\begin{itemize}
\item[(i)]  For all $K\comp X$, $\sup_{p\in K}d_h(u_\eps(p),v_\eps(p)) \to 0$
($\eps\to 0$)
for some (hence every) Riemannian metric $h$ on $Y$.
\item[(ii)] $\forall k\in \N_0\ \forall m\in \N$,
for each chart
   $(V,\vphi)$
   in $X$, each chart $(W,\psi)$ in $Y$, each $L\comp V$
   and each $L'\comp W$:
$$
\sup\limits_{x\in L\cap u_\eps^{-1}(L')\cap v_\eps^{-1}(L')}\!\!\!\!\!\!\!\!\!\!\!\!\!\!\!\!\!\!\!
\|D^{(k)}(\psi\circ u_\eps\circ \vphi^{-1}
- \psi\circ v_\eps\circ \vphi^{-1})(\vphi(p))\|
=O(\eps^m).
$$
\end{itemize}
\et
Due to the c-boundedness condition (i.e., (i) in Definition \ref{modmapmf})
elements $u$ of $\gs[X,Y]$ can model jump discontinuities
but not $\delta$-like singularities; the latter, however, will arise as tangent maps of such $u$
(see below, a list of examples is provided by \cite{gfvm}, 2.8). Inserting a compactly supported point
$\tilde p\in\widetilde X_c$ into $u\in\gs[X,Y]$ yields a well-defined element $[u_\eps(p_\eps)_\eps]
\in\widetilde Y_c$. However, $p\in\widetilde X_c$ characterize $u\in\gs[X,Y]$ only up to equivalence
of order zero ($\sim_0$, to be defined below), i.e., $u(\tilde p)=v(\tilde p)\ \forall \tilde p
\Leftrightarrow (u_\eps)_\eps\sim_0(v_\eps)_\eps$ (Prop.\ 2.14 in \cite{gfvm}).

\bd \label{equiv0def}
We call two elements $(u_\eps)_\eps$, $(v_\eps)_\eps$
of $\esm[X,Y]$ equivalent of order
$0$, $(u_\eps)_\eps \sim_0 (v_\eps)_\eps$ if
they satisfy Definition \ref{equrel} (i) and
(ii) for $k=0$.
\et

We now turn to the definition of generalized vector bundle homomorphisms by
first introducing the respective notions of moderateness and  equivalence.

\bd \label{homgdef}
${\esm}^{\mathrm{VB}}[E,F]$ is the set of all $(u_\eps)_\eps$ 
$\in$ $\mathrm{Hom}(E,F)^I$ satisfying
\begin{itemize}
\item[(i)] $(\underline{u_\eps})_\eps \in \esm[X,Y]$.
\item[(ii)] $\forall k\in \N_0\
\forall (V,\Phi)$
vector bundle chart in $E$,
$\forall (W,\Psi)$ vector bundle chart in $F$,
$\forall L\comp V\
\forall L'\comp W\ \exists N\in \N\ \exists \eps_1>0\
\exists C>0$ with
$$
\|D^{(k)}
(u_{\eps \mathrm{\Psi}\mathrm{\Phi}}^{(2)}(\vphi(p)))\|
\le C\eps^{-N}
$$
for all $\eps<\eps_1$ and all $p\in L\cap\underline{u_\eps}^{-1}(L')$, where $\|\,.\,\|$ 
denotes any matrix norm.
\end{itemize}
\et

\bd \label{homgequ}
$(u_\eps)_\eps$, $(v_\eps)_\eps \in {\esm}^{\mathrm{VB}}[E,F]$
are called $vb$-equivalent, $((u_\eps)_\eps \sim_{vb} (v_\eps)_\eps)$
if
\begin{itemize}
\item[(i)] $(\underline{u_\eps})_\eps \sim (\underline{v_\eps})_\eps$ in
$\esm[X,Y]$.
\item[(ii)] $\forall k\in \N_0\ \forall m\in \N\ \forall (V,\Phi)$
vector bundle chart in
$E$, $\forall (W,\Psi)$ vector bundle chart in $F$,
$\forall L\comp V\ \forall L'\comp W
\ \exists \eps_1>0\ \exists C>0$ such that:
$$
\|D^{(k)}(u_{\eps \mathrm{\Psi}\mathrm{\Phi}}^{(2)}
-v_{\eps \mathrm{\Psi}\mathrm{\Phi}}^{(2)})(\vphi(p))\|
\le C\eps^{m}
$$
for all $\eps<\eps_1$ and all $p\in L\cap\underline{u_\eps}^{-1}(L')
\cap\underline{v_\eps}^{-1}(L')$.
\end{itemize}
\et
We now set $\mathrm{Hom}_{\gs}[E,F] := {\esm}^{\mathrm{VB}}[E,F]\big/\sim_{vb}$.
For $u\in \mathrm{Hom}_{\gs}[E,F]$, $\underline{u} :=[(\underline{u}_\eps)_\eps]$
is a well-defined element of $\gs[X,Y]$ uniquely characterized by $\underline{u}
\circ\pi_X = \pi_Y\circ u$. The tangent map $Tu:=[(Tu_\eps)_\eps]$ of any 
$u\in\gs[X,Y]$ is a well-defined element of $\mathrm{Hom}_{\gs}[TX,TY]$. 

Analogously to the case of manifold valued nets of functions we will need the weaker notion
of $vb$-$0$-equivalence, denoted $\sim_{vb0}$, and defined by
$(u_\eps)_\eps~\sim_{vb0}~(v_\eps)_\eps$, if
$(\underline{u_\eps})_\eps \sim_0 (\underline{v_\eps})_\eps$ and (ii)
above holds for $k=0$.

Finally we define the space $E_c^{\sim_{vb}}$ of compactly bounded generalized
vector bundle points.
\bd \label{vbpoints}
On $(E,X,\pi)$ we define the set of $vb$-moderate generalized points
as consisting of $(e_\eps)_\eps \in E^I$ satisfying
\begin{itemize}
\item[(i)] $\exists K\comp X$ $\exists \eps_0>0$ such that
$(\pi(e_\eps))_\eps \in K$ for all $\eps<\eps_0$.
\item[(ii)] 
$\exists \eps_1>0$ $\exists N\in \N$ $\exists C>0$ such that
$\|e_\eps\| \le C\eps^{-N}$
for all $\eps<\eps_1$ (with the norm again induced by any
Riemannian metric on the base).
\end{itemize}
On this set we introduce the equivalence relation 
$(e_\eps)_\eps \sim_{vb} (e'_\eps)_\eps$ by
\begin{itemize}
\item[(iii)] $(\pi(e_\eps))_\eps \sim (\pi(e'_\eps))_\eps$ in $X^I$
\item[(iv)] $\forall m\in \N 
\ \forall
(W,\Psi)$ vector bundle chart in $E$ $\forall L'\comp W
\ \exists \eps_1>0\ \exists C>0$ such that
$|\bfs{\psi}e_\eps - \bfs{\psi}e'_\eps| \le C\eps^m$
for all $\eps<\eps_1$ 
whenever both $\pi_X(e_\eps)$ and $\pi_X(e'_\eps)$ lie in $L'$.
(Here $\Psi=\big( z_p\mapsto
(\psi(p),\bfs{\psi}(z))\big)$, cf.\ (\ref{vbchart}).)
\end{itemize}
The set of equivalence classes is denoted by $E_c^{\sim_{vb}}$.
\et

Inserting $\tilde p\in\widetilde X_c$ into $u\in\Ga_\gs(X,E)$ yields a well-defined element of
$E_c^{\sim_{vb}}$ and generalized sections are characterized by these point values
(\cite{gfvm}, Th.\ 3.7 (i)). On the other hand $v\in\mathrm{Hom}_\gs[E,F]$ has 
to satisfy an additional condition (equation (6) in \cite{gfvm}) to allow for a well 
defined element $[(v_\eps(e_\eps))_\eps]\in F_c^{\sim_{vb}}$ for $\tilde e\in E_c^{\sim_{vb}}$
(\cite{gfvm}, Prop.\ 3.6 (ii)) and even under this restriction generalized vector bundle homomorphisms 
are characterized by their point values on $E_c^{\sim_{vb}}$ only up to $vb$-$0$-equivalence 
(\cite{gfvm}, Th.\ 3.7 (ii)).

\section{Characterization of generalized functions \\valued in a manifold}\label{emch}

We start this section by establishing simple and global criteria for nets 
$(u_\eps)_\eps\in\CC^\infty(X,Y)^I$ to be moderate. The basic idea is to replace
the chartwise (in $Y$) description of 
Definition \ref{modmapmf} (ii) by composition with smooth (compactly supported) 
functions $f:Y\to\C$. To begin with, we note the following characterization of
the notion of c-boundedness.

\bp\label{propneu}
Let $(u_\eps)_\eps\in\CC^\infty(X,Y)$. The following conditions are equivalent
\begin{itemize}
\item[(i)] $(u_\eps)_\eps$ is c-bounded.
\item[(ii)] $(f\circ u_\eps)_\eps$ is c-bounded for all $f\in \CC^\infty(Y)$.
\item[(iii)] $(f\circ u_\eps)_\eps$ is moderate of order zero for all $f\in \CC^\infty(Y)$, i.e.,
\[
 \forall K\subset\subset X\ \exists N\in\N: \sup_{p\in K}|f\circ u_\eps(p)|=O(\eps^{-N})
\]
for all $f\in \CC^\infty(Y)$.
\item[(iv)] $(u_\eps(x_\eps))_\eps\in Y_c$ for all $(x_\eps)_\eps\in X_c$.
\end{itemize}
\et

\pr (i)$\Rightarrow$(ii)$\Rightarrow$(iii) is clear and so is (i)$\Leftrightarrow$(iv).
We establish (iii) $\Rightarrow$ (i). To begin with suppose $Y$ is non compact and connected.
We then may write $Y=\bigcup_{n\in\N}L_n$ with $L_n$ compact and $L_n$ contained in the interior 
$L_{n+1}^\circ$ of $L_{n+1}$ for all $n$. Suppose (i) is false; 
then $\exists K\comp X\ \forall n\in \N\ \exists \eps_n
\leq 1/n\ \exists p_n\in K$ with $u_{\eps_{n}}(p_n)\not\in L_n$. Without loss of generality we may even
suppose that $u_{\eps_{n}}(p_n)\in L_{n+1}^\circ\setminus L_n$. Now choose $f_n$ in 
$\D(L_{n+1}^\circ\setminus L_n)$ with $f_n(u_{\eps_n}(p_n))=e^{1/\eps_n}$ and let $f:=\sum_{n=1}^\infty
f_n$. Then $f\in\CC^\infty(Y)$ and by (iii) we are guaranteed the existence of $N\in\N$ such that
$|\sup_{p\in K}f(u_\eps(p))|\leq\eps^{-N}$. Hence $|f(u_{\eps_n}(p_n))|\leq\eps_n^{-N}$ but by the above
$e^{1/\eps_n}\leq\eps_n^{-N}$ for large $n$, a contradiction. If $Y$ is not connected we employ a similar
construction taking into account that each connected component of $Y$ is $\sigma$-compact.
\ep\ms

We now have the following characterization of moderateness in $\cinfty(X,Y)^I$.

\bp \label{emchar}
Let $(u_\eps)_\eps \in \cinfty(X,Y)^I$. The following statements are equivalent.
\begin{itemize}
\item[(a)] $(u_\eps)_\eps \in \esm[X,Y]$.
\item[(b)] 
  \begin{itemize}
  \item[(i)] $(u_\eps)_\eps$ is c-bounded.
  \item[(ii)] $(f\circ u_\eps)_\eps \in \esm(X)$ for all $f\in \D(Y)$. 
  \end{itemize}
\item[(c)] 
$(f\circ u_\eps)_\eps \in \esm(X)$ for all $f\in \cinfty(Y)$. 
\end{itemize}
\et

Note that (i) in (b) is necessary. Indeed, let $X=\R=Y$; then  $u_\eps(x)=1/\eps$ 
is not c-bounded but $|f\circ u_\eps(x)|\leq\|f\|_\infty$ for all compactly
supported smooth $f$. 
\ms 

\pr (b)$\Rightarrow$(a): We have to verify (i) and (ii) of 
Definition \ref{modmapmf}.
(i) is identical to (b) (i). Concerning (ii), let $k\in \N$, $(V,\vphi)$, $(W,\psi)$ charts 
in $X$ respectively $Y$ and $L\comp V$, $L'\comp W$. Choose $f\in \D(W)^m$, $f\equiv \psi$ in a
neighborhood of $L'$ and set $f_j:= \pro_j\circ f$. Let $p\in L\cap u_\eps^{-1}(L')$.
Then in a neighborhood (depending on $\eps$) 
of $p$ we have $\psi_j\circ u_\eps = f_j\circ u_\eps$, so
$$
D^{(k)}(\psi_j\circ u_\eps \circ \vphi^{-1})(\vphi(p)) = 
D^{(k)}(f_j\circ u_\eps \circ \vphi^{-1})(\vphi(p))\,,
$$ 
from which the claim follows.

\noindent (a)$\Rightarrow$(c): We have to deduce (c) from 
Definition \ref{modmapmf}. Let $f\in
\CC^\infty(Y)$ and $K\comp X$. We may assume without loss of generality that $K\comp V$ for some
chart $(V,\vphi)$ in $X$. Since $(u_\eps)_\eps$ is c-bounded we may choose $K'\comp Y$ and
$\eps_0>0$ such that $u_\eps(K)\subseteq K'$ for all $\eps<\eps_0$. Cover $K'$ by 
charts $(W_l,\psi_l)$ in $Y$ ($1\le l\le s$) and write $K'=\bigcup_{l=1}^s K_l'$
with $K_l'\comp W_l$ for each $l$.
Let $p\in K$ and choose $l\in \{1,\dots,s\}$ such that $u_\eps(p)\in K_l'$.
Then $f\circ u_\eps = (f\circ \psi_l^{-1})\circ (\psi_l\circ u_\eps)$ in a neighborhood
of $p$ (depending on $\eps$). Applying 
Definition \ref{modmapmf} (ii) to $L:=K$, $(V,\vphi)$, $L':=K_l'$, $(W_l,\psi_l)$ 
we obtain $\eps_1 = \min_{1\le l\le s}\eps_1^l$, 
$N=\max_{1\le l\le s}N_l$, and $C=\max_{1\le l\le s}C_l$ such that
$$
\sup_{p\in K} \|D^{(k)}(f\circ u_\eps \circ\vphi^{-1})(\vphi(p))\| \le C \eps^{-N}
$$ 
for $\eps<\eps_1$, so indeed $(f\circ u_\eps)_\eps \in \esm(X)$.

\noindent (c)$\Rightarrow$(b) is immediate from Proposition\ref{propneu}.
\ep\ms

A similar characterization can be derived for the equivalence relation $\sim$
on $\esm[X,Y]$ (see (c), (d) in Theorem \ref{equivchar} below). 
Moreover, this characterization provides an affirmative answer to
the question raised in \cite{gfvm}, Remark 2.11 whether a characterization of equivalence
of elements of spaces of Colombeau generalized functions which does not resort to
derivatives of representatives (as established for practically all (scalar) variants of the
construction in \cite{found}, Th.\ 13.1) was also attainable in the context of
manifold-valued generalized functions.

\bt \label{equivchar} 
Let $(u_\eps)_\eps$, $(v_\eps)_\eps \in \esm[X,Y]$. The following statements are
equivalent:
\begin{itemize}
\item[(i)] $(u_\eps)_\eps \sim (v_\eps)_\eps$.
\item[(ii)] $(u_\eps)_\eps \sim_0 (v_\eps)_\eps$.
\item[(iii)] For every
Riemannian metric $h$ on $Y$, every $m\in \N$ and every $K\comp X$,  
$$
\sup_{p\in K}d_h(u_\eps(p),v_\eps(p)) = O(\eps^m) \qquad (\eps\to 0)\,.
$$ 
\item[(iv)] $(f\circ u_\eps - f\circ v_\eps)_\eps \in \ns(X)$ for all $f\in \D(Y)$.
\item[(v)] $(f\circ u_\eps - f\circ v_\eps)_\eps \in \ns(X)$ for all $f\in \cinfty(Y)$.
\end{itemize}
\et

\pr Clearly, (i)$\Rightarrow$(ii).

\noindent (ii)$\Leftrightarrow$(iii): see \cite{gfvm}, Th.\ 2.10.

\noindent (ii)$\Rightarrow$(iv): Let $K\comp X$, $f\in \D(Y)$. 
By \cite{ndg}, (3) it suffices to show that 
$$
\sup_{p\in K}|f\circ u_\eps(p) -f\circ v_\eps(p)| = O(\eps^m)
$$ 
for each $m\in \N$. To see this, choose $K'\comp Y$ and $\eps_0>0$ such that 
$u_\eps(K)\cup v_\eps(K) \subseteq K'$ for $\eps<\eps_0$. Cover $K'$ by open sets
$W_l'$, $\overline{W_l'}\comp W_l$ for charts $(W_l,\psi_l)$ 
($1\le l\le s$). 
From (ii) it follows that for each fixed sufficiently small $\eps$ and each $p\in K$ there exists 
$l\in \{1,\dots,s\}$ such that both $u_\eps(p)$ and $v_\eps(p)$ are contained
in $W_l'$. In this case we have
\beas
|(f\circ u_\eps-f\circ v_\eps)(p)| &=& |((f\circ \psi_l^{-1})\circ(\psi_l\circ u_\eps)
-(f\circ \psi_l^{-1})\circ(\psi_l\circ v_\eps))(p)|\\ 
&\le& C |(\psi_l\circ v_\eps - \psi_l\circ u_\eps)(p)| 
\eeas
by \cite{gfvm}, Lemma 2.5, where $C$ is independent of $\eps$. Hence the claim follows
from 
Definition \ref{equiv0def}.

\noindent (iv)$\Rightarrow$(i): We first have to show that for every Riemannian metric
$h$ on $Y$ and any $K\comp X$ we have $\sup_{p\in K}d_h(u_\eps(p),v_\eps(p))\to 0$ as 
$\eps\to 0$. Suppose to the contrary that 
\begin{equation}\label{indirect}
    \exists K\comp X \ \exists \delta>0 \ \forall k\in \N \ \exists \eps_k<\frac{1}{k}
\ \exists p_k\in K: \quad 
d_h(u_{\eps_k}(p_k),v_{\eps_k}(p_k)) \ge \delta\,.
\end{equation}
Since $(u_\eps)_\eps$ and $(v_\eps)_\eps$ are c-bounded, there exists $K'\comp Y$ such that
$u_\eps(K)\cup v_\eps(K)\subseteq K'$ for $\eps$ small. Thus, without loss of generality
we may suppose that $u_{\eps_k}(p_k)\to q_1$ and $v_{\eps_k}(p_k)\to q_2$ with $q_i\in K'$
for $i=1,2$. By (\ref{indirect}), $q_1\not= q_2$. Choose $f\in \D(Y)$ with $f(q_1)=1$, 
$f(q_2)=0$. But then
$$
1 = \lim_{k\to\infty}|f(u_{\eps_k}(p_k))-f(v_{\eps_k}(p_k))| \le \sup_{p\in K}
|f(u_{\eps_k}(p))-f(v_{\eps_k}(p))| = O(\eps_k^m)
$$
for all $m$, a contradiction. 

Finally, we have to establish property (ii) of 
Definition \ref{equrel}. Thus let
$L\comp V$, $(V,\vphi)$ a chart in $X$, $L'\comp W$, $(W,\psi)$ a chart in $Y$.
Fix $j\in \{1,\dots,m\}$ and choose $f_j\in \D(W)$ such that $f_j\equiv\psi_j$
in a neighborhood of $L'$. Then for any $p\in L\cap u_\eps^{-1}(L')\cap v_\eps^{-1}(L')$
there exists a neighborhood (depending on $\eps$) on which $\psi_j\circ u_\eps = f_j\circ u_\eps$ 
and $\psi_j\circ v_\eps = f_j\circ v_\eps$. Hence,
since $(f_j\circ u_\eps - f_j\circ v_\eps)_\eps \in \ns(X)$,  
for each $k$ in $\N$ we have
\beas
&& \|D^{(k)}(\psi_j\circ u_\eps\circ \vphi^{-1}
- \psi_j\circ v_\eps\circ \vphi^{-1})(\vphi(p))\| \\ 
&& \hphantom{mmmm} \le 
\sup_{p'\in L}\|D^{(k)}(f_j\circ u_\eps\circ \vphi^{-1}
- f_j\circ v_\eps\circ \vphi^{-1})(\vphi(p'))\| = O(\eps^m)
\eeas
for each $m$, yielding the claim.

\noindent (iv)$\Leftrightarrow$(v) follows as in Proposition \ref{emchar}.
\ep\ms

The above characterization result has several important consequences. 
First, we obtain a characterization of equivalence of compactly supported generalized points.
\bc Let $(p_\eps)_\eps$, $(q_\eps)_\eps \in X^I$ be compactly supported. The following
statements are equivalent:
\begin{itemize}
\item[(i)] $(p_\eps)_\eps \sim (q_\eps)_\eps$.
\item[(ii)] $(f(p_\eps))_\eps \sim (f(q_\eps))_\eps$ for all 
$f\in \D(X)$ (respectively all $f\in \cinfty(X)$). 
\item[(iii)] $|f(p_\eps)-f(q_\eps)| = O(\eps^m)$ as $\eps\to 0$ for all $m\in \N$
and all $f\in \D(X)$ (respectively all $f\in \cinfty(X)$).
\end{itemize}
\et

\pr (i)$\Leftrightarrow$(iii): 
Setting $u_\eps(p)\equiv p_\eps$ and $v_\eps(p)\equiv q_\eps$, the nets $(u_\eps)_\eps$,
$(v_\eps)_\eps$ are contained in $\esm[X,X]$. It is clear from the definitions that
$(u_\eps)_\eps \sim (v_\eps)_\eps$ in $\esm[X,X]$ if and only if $(p_\eps)_\eps \sim
(q_\eps)_\eps$. Thus an application of Theorem \ref{equivchar}, (i)$\Leftrightarrow$(iv) 
gives the result.

\noindent(ii)$\Leftrightarrow$(iii) is immediate from the definition of equivalence
of generalized numbers. 
\ep\ms

Based on Theorem \ref{equivchar}, the following result provides a point value 
characterization of elements of $\gs[X,Y]$. 

\bt \label{pvchar} Let $u,\, v \in \gs[X,Y]$. Then $u=v$ if and only if $u(\tilde p) = 
v(\tilde p)$ for all $\tilde p \in \widetilde X_c$.  
\et

\pr Immediate from (a)$\Leftrightarrow$(b) in Theorem \ref{equivchar} and 
\cite{gfvm}, Prop.\ 2.14.
\ep\ms

Moreover, Theorem \ref{equivchar} can be used to show that composition of manifold-valued 
generalized functions can be carried out unrestrictedly (i.e., without any additional
assumptions, as, e.g., condition (6) in \cite{gfvm}).
\bt \label{generalcomp} Let 
$u=[(u_\eps)_\eps] \in \gs[X,Y]$, $v=[(v_\eps)_\eps] 
\in \gs[Y,Z]$. Then $v\circ u := [(v_\eps\circ u_\eps)_\eps]$ is a well-defined
element of $\gs[X,Z]$. 
\et

\pr The proof of \cite{gfvm}, Th.\ 2.16 shows that $(v_\eps\circ u_\eps)_\eps 
\in \esm[X,Y]$. Moreover, it also establishes that $(u_\eps)_\eps \sim 
(u_\eps')_\eps$ implies $(v_\eps\circ u_\eps)\sim_0 (v_\eps\circ u_\eps')$ and
that $(v_\eps)_\eps \sim (v_\eps')_\eps$ implies $(v_\eps\circ u_\eps)\sim_0 (v_\eps'
\circ u_\eps)$. Hence by Theorem \ref{equivchar}, the class of $(v_\eps\circ u_\eps)_\eps$
is well-defined.
\ep\ms

\section{Characterization of generalized vector bundle homomorphisms}\label{homch}

The aim of this section is to derive characterizations similar to Proposition \ref{emchar} and
Theorem \ref{equivchar} also for generalized vector bundle homomorphisms. 
In the following considerations for any  vector bundle $E$ over $X$ 
the role of the  ``model spaces'' $\esm(X)$ respectively $\ns(X)$ (cf.\ Proposition 
\ref{emchar} (b), (c) and Theorem \ref{equivchar} (iv), (v)) will be played
by $\esm^{VB}[E,\R\times \R^{m'}]$ respectively $\sim_{vb}$ on that space.
In order to allow a smooth presentation of the following results let us therefore
examine the simplifications that ensue from the particularly simple structure
of the range space $\R\times \R^{m'}$. Let $(u_\eps)_\eps \in 
\mathrm{Hom}(E,\R\times\R^{m'})$ and let $(V,\Phi)$ be a vector bundle chart in
$E$. Then using the trivial vector bundle chart $\mathrm{id}$ on $\R\times\R^{m'}$
we may write
\begin{equation} \label{locformfu}
u_\eps\circ\Phi^{-1}(x,\xi) = 
(u_{\eps \mathrm{id}\Phi}^{(1)}(x),
u_{\eps\mathrm{id}\Phi}^{(2)}(x)\cdot\xi)\,.
\end{equation}
It follows that $(u_\eps)_\eps \in \esm^{VB}[E,\R\times\R^{m'}]$ if and only if
$(\underline{u_\eps})_\eps$ is c-bounded and for each vector bundle chart $(V,\Phi)$,
$u_{\eps_\mathrm{id}\Phi}^{(1)} \in \esm(\vphi(V))$ and $u_{\eps_\mathrm{id}\Phi}^{(2)}
\in \esm(\vphi(V))^{(m')^2}$. Moreover, it follows that for 
$(u_\eps)_\eps,\, (v_\eps)_\eps \in \esm^{VB}[E,\R\times\R^{m'}]$, $(u_\eps)_\eps
\sim_{vb} (v_\eps)_\eps$ is equivalent to $u_{\eps_\mathrm{id}\Phi}^{(1)}
- v_{\eps_\mathrm{id}\Phi}^{(1)} \in \ns(\vphi(V))$ and 
$u_{\eps_\mathrm{id}\Phi}^{(2)} - v_{\eps_\mathrm{id}\Phi}^{(2)} \in \ns(\vphi(V))$ for
each vector bundle chart $(V,\Phi)$.
By \cite{found}, Th.\ 13.1 and the remark following it, this in turn is equivalent
to $(u_\eps)_\eps \sim_{vb0} (v_\eps)_\eps$, so $\sim_{vb}$ and $\sim_{vb0}$ are
the same for the model space $\esm^{VB}[E,\R\times\R^{m'}]$. We will extend
the validity of this observation to general $\esm[E,F]$ in Theorem \ref{equivvbchar}
below.

For a vector bundle $E \to X$ we denote by $\mathrm{Hom}_c(E,\R\times\R^{m'})$ the 
set of all vector bundle homomorphisms $f$ such that $\underline f:X\to\R$ has
compact support.
\bp \label{emvbchar} 
Let $(u_\eps)_\eps \in \mathrm{Hom}(E,F)^I$. The following statements are
equivalent:
\begin{itemize}
\item[(a)] $(u_\eps)_\eps \in \esm^{VB}[E,F]$.
\item[(b)] 
  \begin{itemize}
  \item[(i)] $(\underline{u_\eps})_\eps$ is c-bounded.
  \item[(ii)] $(\hat f \circ u_\eps)_\eps \in \esm^{VB}(E,\R\times\R^{m'})$ 
              for all $\hat f \in \mathrm{Hom}_c(F,\R\times \R^{m'})$.
  \end{itemize}
\item[(c)] 
$(\hat f \circ u_\eps)_\eps \in \esm^{VB}(E,\R\times\R^{m'})$ 
              for all $\hat f \in \mathrm{Hom}(F,\R\times \R^{m'})$.
\end{itemize}
\et
\pr (a)$\Rightarrow$(c): 
By (\ref{locformfu})
we have to show that for each $\hat f \in \mathrm{Hom}(F,\R\times \R^{m'})$
and each vector bundle chart $(V,\Phi)$ of $E$ we have 
$((\hat f\circ u_\eps)_{\mathrm{id}\Phi}^{(1)})_\eps\in
\esm(\vphi(V))$ and $((\hat f\circ u_\eps)_{\mathrm{id}\Phi}^{(2)})_\eps \in
\esm(\vphi(V))^{(m')^2}$. Due to $(\hat f\circ u_\eps)_{\mathrm{id}\Phi}^{(1)}
= \underline{\hat f}\circ \underline{u_\eps}\circ\vphi^{-1}$, the first of these
claims follows from Proposition \ref{emchar} since $(\underline{u_\eps})_\eps\in 
\esm[X,Y]$ and $\underline{\hat f} \in \cinfty(Y)$. 
Concerning the second one, let
$L\comp V$ and choose $\eps_0>0$, $L'\comp Y$ such that $\underline{u_\eps}(L)
\subseteq L'$ for $\eps<\eps_0$. Cover $L'$ by 
vector bundle charts $(W_l,\Psi_l)$ of $F$ ($1\le l\le s$) and write 
$L=\bigcup_{l=1}^s L_l'$ with $L_l'\comp W_l$. 
Choose $\eps_1<\eps_0$ such that property (ii) of 
Definition \ref{homgdef} is satisfied for $(L,V,\Phi)$ and $(L_1',W_1,\Psi_1)$,
\dots, $(L_s',W_s,\Psi_s)$ simultaneously. Let $\eps<\eps_1$ and $p\in L$ with
$u_\eps(p)\in L_l'$. Then for $p'$ in a neighborhood (depending on $\eps$) of $p$ we have
\begin{equation} \label{locformfu2}
(\hat f\circ u_\eps)_{\mathrm{id}\Phi}^{(2)}(\vphi(p'))
= (\hat f \circ \Psi_l^{-1})^{(2)}(\psi_l\circ\underline{u_\eps}(p'))\cdot 
u_{\eps\Psi_l\Phi}^{(2)}(\vphi(p'))\,,
\end{equation}
so the claim follows.

\noindent(c)$\Rightarrow$(b): (i) follows from Proposition \ref{propneu} while (ii) is immediate.

\noindent(b)$\Rightarrow$(a) We have to establish properties (i) and (ii) of 
Definition \ref{homgdef}. Since $(\hat f\circ u_\eps)_{\mathrm{id}\Phi}^{(1)}
= \underline{\hat f}\circ \underline{u_\eps}\circ\vphi^{-1}$ for any $\hat f\in
\mathrm{Hom}(F,\R\times\R^{m'})$, (i) follows from Proposition \ref{emchar}.
Let $L\comp V$, $(V,\Phi)$ a vector bundle chart in $E$,
$L'\comp W$, $(W,\Psi)$ a vector bundle chart in $F$ and choose an open
neighborhood $U$ of $L'$ whose compact closure is contained in $W$. Choose any
$l\in \{1,\dots, m\}$ and let $\hat f_l\in \mathrm{Hom}_c(F,\R\times\R^{m'})$
such that 
$$
\hat f_l|_{\pi_Y^{-1}(U)} = (\pro_l \times \mathrm{id}_{\R^{m'}}) 
\circ \Psi|_{\pi_Y^{-1}(U)}
$$
where $\pro_l:\R^m\to \R$. Let $p\in L\cap \underline{u_\eps}^{-1}(L')$. Then for $p'$ 
in a suitable neighborhood of $p$ we have
\begin{equation}\label{dingsda}
(\hat f_l
\circ u_\eps\circ\Phi^{-1})^{(2)}(\vphi(p')) = 
u_{\eps\Psi\Phi}^{(2)}(\vphi(p'))\,,
\end{equation}
so the desired estimates of $D^{(k)}(u_{\eps\Psi\Phi}^{(2)})(\vphi(p))$ follow,
thereby establishing property (ii) of 
Definition \ref{homgdef}.
\ep\ms

We now characterize the equivalence relation $\sim_{vb}$ on $\esm^{VB}[E,F]$,
simultaneously establishing that vb-equivalence is in fact identical to 
vb-$0$-equivalence.
\bt \label{equivvbchar}
Let $(u_\eps)_\eps$, $(v_\eps)_\eps \in \esm^{VB}(E,F)$. The following statements are 
equivalent:
\begin{itemize}
\item[(i)] $(u_\eps)_\eps \sim_{vb} (v_\eps)_\eps$.
\item[(ii)] $(u_\eps)_\eps \sim_{vb0} (v_\eps)_\eps$.
\item[(iii)] $(\hat f\circ u_\eps)\sim_{vb} (\hat f\circ v_\eps)_\eps$ in 
$\esm^{VB}[E,\R\times \R^{m'}]$ for all $\hat f\in \mathrm{Hom}_c(E,\R\times \R^{m'})$. 
\item[(iv)] $(\hat f\circ u_\eps)\sim_{vb} (\hat f\circ v_\eps)_\eps$ in 
$\esm^{VB}[E,\R\times \R^{m'}]$ for all $\hat f\in \mathrm{Hom}(E,\R\times \R^{m'})$. 
\end{itemize}
\et
\pr (i)$\Rightarrow$(ii) is clear.

\noindent(ii)$\Rightarrow$(iv): Let $\hat f\in \mathrm{Hom}(E,\R\times \R^{m'})$
and $L\comp V$, $(V,\Phi)$ a vector bundle chart in $E$. To prove (iv), by 
the remarks following (\ref{locformfu}) it suffices to establish the $\ns$-estimates of 
order $0$ on $L$ for both $(\hat f\circ u_\eps-\hat f \circ v_\eps)_{\mathrm{id}\Phi}^{(1)}$ and 
$(\hat f\circ u_\eps-\hat f \circ v_\eps)_{\mathrm{id}\Phi}^{(2)}$. 
Choose $K'\comp Y$ and $\eps_1>0$ such that 
$\underline{u_\eps}(L)\cup \underline{v_\eps}(L) \subseteq K'$ for 
$\eps<\eps_1$. Cover $K'$ by open sets
$W_l'$, $\overline{W_l'}\comp W_l$ ($1\le l\le s$), where $(W_l,\Psi_l)$
are vector bundle charts in $F$. 
Since for any Riemannian metric $h$ on $Y$, 
$\sup_{p\in K}d_h(\underline{u_\eps}(p),\underline{v_\eps}(p))\to 0$ as $\eps\to 0$, 
for each $\eps$ small and each $p\in K$ there exists 
$l\in \{1,\dots,s\}$ such that both $\underline{u_\eps}(p)$ and $\underline{v_\eps}(p)$ 
are contained in $W_l'$. Now the estimate for 
\begin{equation} \label{eq4}
(\hat f\circ u_\eps-\hat f \circ v_\eps)_{\mathrm{id}\Phi}^{(1)}(\vphi(p)) =
(\underline{\hat f}\circ \underline{u_\eps} - 
\underline{\hat f}\circ \underline{v_\eps})\circ\vphi^{-1}(p)
\end{equation} 
follows from Theorem \ref{equivchar} (ii) $\Rightarrow$ (v), while that for 
$(\hat f\circ u_\eps-\hat f \circ v_\eps)_{\mathrm{id}\Phi}^{(2)}(\vphi(p))$
is derived from a representation as in (\ref{locformfu2}) using 
Definition \ref{homgequ} (ii) (for $k=0$) and Lemma 2.5 of \cite{gfvm}: again the constants appearing due to the
application of that result can be chosen independently of $\eps$ since the 
domain of $(\hat f \circ \Psi_l^{-1})^{(2)}$
is an open neighborhood of $\psi_l(\overline{W_l'})$.

\noindent(iv)$\Rightarrow$(iii): Clear.

\noindent(iii)$\Rightarrow$(i) Using (\ref{eq4}), property (i) of 
Definition \ref{homgequ} follows from Theorem \ref{equivchar} (v) $\Rightarrow$ (i). 
Finally, property (ii) of that
definition is established by employing an $\hat f_l$ as in the proof of
Proposition \ref{emvbchar}, (b)$\Rightarrow$(a), using a representation of both 
$u_{\eps\Psi\Phi}^{(2)}$ and $v_{\eps\Psi\Phi}^{(2)}$ as in (\ref{dingsda}).
\ep\ms

The above result can in turn be used to derive the following characterization result
for generalized vector bundle points.
\bc \label{equivvbcharpoint}
Let $(e_\eps)_\eps$, $(e_\eps')_\eps \in E^I$ be vb-moderate. 
The following statements are equivalent:
\begin{itemize}
\item[(i)] $(e_\eps)_\eps \sim_{vb} (e_\eps')_\eps$.
\item[(ii)] $(\hat f(e_\eps))_\eps \sim_{vb} (\hat f(e_\eps'))_\eps$ in
$\R \times \R^{n'}$ for each $\hat f\in 
\mathrm{Hom}_c(E,\R\times \R^{n'})$ (respectively $\in \mathrm{Hom}(E,\R\times \R^{n'})$). 
\item[(iii)] $|\hat f(e_\eps) -\hat f(e_\eps')| =O(\eps^m)$ as $\eps\to 0$ for all
$m\in \N$ and all $\hat f\in \mathrm{Hom}_c(E,\R\times \R^{n'})$
(respectively $\in \mathrm{Hom}(E,\R\times \R^{n'})$). 
\end{itemize}
\et

\pr (i)$\Leftrightarrow$(ii) follows from applying Theorem \ref{equivvbchar} to 
$u_\eps\equiv e_\eps$, $v_\eps\equiv e_\eps'$.

\noindent(ii)$\Leftrightarrow$(iii) is immediate from the definitions.
\ep\ms

We can employ this result in order to get rid of the additional technical assumption 
in the statement of \cite{gfvm}, Prop.\ 3.6 (ii) (called (6) there). Hence 
generalized vector bundle points may be inserted into generalized vector bundle homomorphisms
unrestrictedly. 
\bc \label{pointinsert}
Let $u= [(u_\eps)_\eps] \in \mathrm{Hom}_{\gs}[E,F]$, $\tilde e = [(e_\eps)_\eps] \in 
E_c^{\sim_{vb}}$. Then $u(\tilde e) := [(u_\eps(e_\eps))_\eps]$ is a well-defined
element of of $F_c^{\sim_{vb}}$.
\et
\pr The proof of \cite{gfvm}, Prop.\ 3.6 (ii) shows (even under the more general
assumptions made here) that $(u_\eps(e_\eps))_\eps$ is vb-moderate and that 
$(u_\eps)_\eps\sim_{vb}(u_\eps')_\eps$ implies $(u_\eps(e_\eps))_\eps \sim_{vb}
(u_\eps'(e_\eps))_\eps$ for any vb-moderate $(e_\eps)_\eps$. To finish the proof
it remains to show that $(e_\eps)_\eps \sim (e_\eps')_\eps$ implies $(u_\eps(e_\eps))_\eps
\sim_{vb} (u_\eps(e_\eps'))_\eps$ even for general $(u_\eps)_\eps \in \esm^{VB}[E,F]$.
Let $\hat f\in \mathrm{Hom}_c(F,\R\times \R^{m'})$. By Corollary \ref{equivvbcharpoint}
(ii) it suffices to establish $((\hat f\circ u_\eps)(e_\eps))_\eps \sim_{vb}
((\hat f\circ u_\eps)(e_\eps'))_\eps$. Now $(\hat f\circ u_\eps)_\eps \in \esm^{VB}[E,\R\times 
\R^{m'}]$ by Proposition \ref{emvbchar} and by Ex.\ 2.15 (i) of \cite{gfvm}, 
$(\underline{\hat f\circ u_\eps})_\eps$ satisfies (6) from \cite{gfvm}. 
Thus the claim follows from \cite{gfvm}, Prop.\ 3.6 (ii).
\ep\ms

The above result enables us to conclude the following point value characterization
of generalized vector bundle homomorphisms.
\bt \label{pvcharhom} 
Let $u,\,  v \in \mathrm{Hom}_{\gs}[E,F]$. Then
$u=v$ if and only if $u(\tilde e) = v(\tilde e)$ in $F_c^{\sim_{vb}}$ for all
$\tilde e \in E_c^{\sim_{vb}}$.
\et

\pr By Corollary \ref{pointinsert}, the point values of $u$ and $v$ at $\tilde e$
are well-defined. In particular, $u=v$ entails $u(\tilde e) = v(\tilde e)$ for
all $\tilde e \in E_c^{\sim_{vb}}$. The converse direction follows from the proof
of \cite{gfvm}, Th.\ 3.7 (ii) and an application of Theorem \ref{equivvbchar},
(i)$\Leftrightarrow$(ii).
\ep\ms

Furthermore, it now follows that composition of generalized vector bundle 
homomorphisms can be carried out unrestrictedly:
\bt \label{pvvbth}
Let $u=[(u_\eps)_\eps] \in \mathrm{Hom}_{\gs}[E,F]$, $v=[(v_\eps)_\eps] 
\in \mathrm{Hom}_{\gs}[F,G]$. Then $v\circ u := [(v_\eps\circ u_\eps)_\eps]$ is
a well-defined element of $\mathrm{Hom}_{\gs}[E,G]$.  
\et

\pr In order to show that $(v_\eps\circ u_\eps)_\eps \in \esm^{VB}[E,G]$, by 
Proposition \ref{emvbchar} it suffices to show that $(\hat f\circ v_\eps\circ u_\eps)
\in \esm^{VB}[E,\R\times \R^{k'}]$ for each $\hat f \in 
\mathrm{Hom}_c[G,\R\times\R^{k'}]$. Moreover, also by Proposition \ref{emvbchar}
we know that for each $\hat f\in\mathrm{Hom}_{c}[G,\R\times\R^{k'}]$ 
$(\hat f \circ v_\eps)_\eps$ is an element of $\esm^{VB}[F,\R\times \R^{k'}]$. 
We may therefore assume without loss of generality that $G=\R\times\R^{k'}$. 
$(\underline{v_\eps\circ u_\eps})_\eps = (\underline{v_\eps}\circ \underline{u_\eps})_\eps
\in \esm[X,Z]$ by Theorem \ref{generalcomp}. Let $(V,\Phi)$ be any vector bundle chart 
in $E$. By the remarks following (\ref{locformfu}) it remains to establish the 
$\esm$-estimates for $(v_\eps\circ u_\eps)_{\mathrm{id}\Phi}^{(2)}\circ \vphi$ 
on every $K\comp V$.
Choose $\eps_1>0$ and $K'\comp Y$ such that $\underline{u_\eps}(K)\subseteq K'$ for all
$\eps<\eps_1$. Cover $K'$ by open sets $W_1$, \dots, $W_l$ 
for vector bundle charts $(W_i,\Psi_i)$ and write $K'=\bigcup_{i=1}^l K_i'$ 
with $K_i'\comp W_i$ ($1\le i\le l$). Let $p \in K$ and $\eps<\eps_1$. Choose $i\in
\{1,\dots,l\}$ such that $\underline{u_\eps}(p) \in K_i'$. Then for $p'$ in a 
neighborhood of $p$ we have  
$$
(v_\eps\circ u_\eps)_{\mathrm{id}\Phi}^{(2)}(\vphi(p')) =
v_{\eps\mathrm{id}\Psi_i}^{(2)}(\psi_i\circ \underline{u_\eps}(p'))\cdot
u_{\eps\Psi_i\Phi}^{(2)}(\vphi(p'))\,.  
$$
Hence the desired estimates for $D^{(k)}((v_\eps\circ u_\eps)_{\mathrm{id}\Phi}^{(2)})
(\vphi(p))$ follow from the moderate\-ness-estimates of $v_{\eps\mathrm{id}\Psi_i}^{(2)}$
on $\psi_i(K_i')$ and property (ii) of 
Definition \ref{homgdef}, applied to 
$u_{\eps\Psi_i\Phi}^{(2)}$ and the set of data $(K,V,\Phi)$, $(K_i',W_i,\Psi_i)$.

To show that $[(v_\eps\circ u_\eps)_\eps]$ is well-defined we resort to the point value
characterization derived in Theorem \ref{pvvbth}: if $(u_\eps)_\eps \sim_{vb} (u_\eps')_\eps$
then for any $\tilde e =[(e_\eps)_\eps]$ we have $(u_\eps(e_\eps))_\eps\sim_{vb}
(u_\eps'(e_\eps))_\eps$, hence $(v_\eps\circ u_\eps(e_\eps))_\eps\sim_{vb}
(v_\eps\circ u_\eps(e_\eps))_\eps$ and therefore $(v_\eps\circ u_\eps)_\eps\sim_{vb}
(v_\eps\circ u_\eps)_\eps$. Similarly, $(v_\eps)_\eps\sim_{vb} (v_\eps')_\eps$ entails
$(v_\eps\circ u_\eps)_\eps\sim_{vb} (v_\eps\circ u_\eps)_\eps$, which concludes
the proof.
\ep\ms

\section{Generalized vector bundle homomorphisms on a generalized mapping}\label{hybch}

This section is devoted to an analysis of the hybrid Colombeau spaces introduced in \cite{curvature},
Sec.\ 4. For the convenience of the reader we recall the basic definitions, namely of
the space of moderate nets in ${\cal E}^h[X,F]=\CC^\infty(X,F)^I$ and the appropriate notion 
of equivalence in this space.

\bd \label{hybridmod} 
We define $\esmh[X,F]$ as the set of all nets 
$(u_\eps)_\eps \in{\cal E}^h[X,F]$ satisfying 
\begin{itemize}
\item[(i)] 
$\forall K\comp X\ \exists K'\comp Y\ \exists \eps_0>0\ \forall \eps<\eps_0$
$\underline{u_\eps}(K) \subseteq K'$.
\item[(ii)] $\forall k\in \N_0
\ \forall (V,\vphi)$ chart in $X$ $\forall (W,\Psi)$ vector bundle
chart in $F$ $\forall L\comp V\ \forall L'\comp W
\ \exists N\in \N\ \exists \eps_1>0\ \exists C>0$ such that
$$\|D^{(k)}(\Psi\circ u_\eps \circ \vphi^{-1})(\vphi(p))\| \le
C\eps^{-N}$$
for each $\eps<\eps_1$ and each $p\in
L\cap\underline{u_\eps}^{-1}(L')$.
\end{itemize}
\et

Note that in particular, $(u_\eps)_\eps \in \esmh[X,F]$ implies $(\underline{u_\eps})_\eps
\in \esm[X,Y]$.

\bd \label{hybridequiv} 
$(u_\eps)_\eps$, $(v_\eps)_\eps \in \esmh[X,F]$ are called equivalent,
$(u_\eps)_\eps \simh (v_\eps)_\eps$, if the following conditions
are satisfied:
\begin{itemize}
\item[(i)] For each $K\comp X$,
$\sup_{p\in X}d_h(\underline{u_\eps},\underline{v_\eps}) \to 0$ 
\item[(ii)] $\forall k\in \N_0\ \forall m\in \N\ \forall (V,\vphi)$
chart in $X$, $\forall (W,\Psi)$ vector bundle
chart in $F$, $\forall L\comp V\ \forall L'\comp W
\ \exists \eps_1>0\ \exists C>0$ such that
$$\|D^{(k)}(\Psi\circ u_\eps\circ\vphi^{-1} - \Psi\circ v_\eps\circ\vphi^{-1})
(\vphi(p))\| \le C\eps^m$$
for each $\eps<\eps_1$ and each $p\in
L\cap\underline{u_\eps}^{-1}(L')\cap
\underline{v_\eps}^{-1}(L')$.
\end{itemize}
\et

The space of hybrid Colombeau generalized functions from the manifold $X$ 
into the vector bundle $F$ as usual is defined by $\gsh[X,F] := \esmh[X,F]\big/\simh$.
We now start by establishing characterizations of moderateness respectively equivalence 
analogous to the previous cases.


\bp \label{emhybridchar} 
Let $(u_\eps)_\eps \in \cinfty(X,F)^I$. The following statements are equivalent:
\begin{itemize}
\item[(a)] $(u_\eps)_\eps \in \esm^h[X,F]$.
\item[(b)] 
  \begin{itemize}
  \item[(i)] $(\underline{u_\eps})_\eps$ is c-bounded. 
  \item[(ii)] $(\hat f \circ u_\eps)_\eps \in \esm^h[X,\R\times\R^{m'}]$
              for all $\hat f \in \mathrm{Hom}_c(F,\R\times\R^{m'})$. 
  \end{itemize}
\item[(c)] 
$(\hat f \circ u_\eps)_\eps \in \esm^h[X,\R\times\R^{m'}]$
              for all $\hat f \in \mathrm{Hom}(F,\R\times\R^{m'})$. 
\end{itemize}
\et

\pr
\noindent (a)$\Rightarrow$(c): 
We have to show that for each $\hat f = (\underline{\hat f}\circ \pi_Y,\hat f_2) \in 
\mathrm{Hom}(F,\R\times \R^{m'})$ (with $\hat f_2 = \pro_2\circ \hat f$, $\pro_2:\R\times 
\R^{m'}\to \R^{m'}$) and each chart $(V,\vphi)$ in $X$ we have 
$\underline{\hat f} \circ \underline{u_\eps}\circ \vphi^{-1}\in
\esm(\vphi(V))$ and $(\hat f_2 \circ u_\eps \circ \vphi^{-1})_\eps \in
\esm(\vphi(V))^{m'}$. The first 
claim follows from Proposition \ref{emchar} since $(\underline{u_\eps})_\eps\in 
\esm[X,Y]$ and $\underline{\hat f} \in \cinfty(Y)$. Let
$L\comp V$ and choose $\eps_0>0$, $L'\comp Y$ such that $\underline{u_\eps}(L)
\subseteq L'$ for $\eps<\eps_0$. Choose a covering of $L'$ by 
vector bundle charts $(W_l,\Psi_l)$ of $F$ ($1\le l\le s$) and write 
$L=\bigcup_{l=1}^s L_l'$ with $L_l'\comp W_l$. There exists
$\eps_1<\eps_0$ such that property (ii) of 
Definition \ref{hybridmod} is satisfied for $(L,V,\vphi)$ and $(L_1',W_1,\Psi_1)$,
\dots, $(L_s',W_s,\Psi_s)$ simultaneously. Let $\eps<\eps_1$ and $p\in L$ with
$u_\eps(p)\in L_l'$. Then 
in a neighborhood (depending on $\eps$) of 
$\vphi(p)$ we have
\begin{equation} \label{locformfu2new}
\hat f_2 \circ u_\eps \circ \vphi^{-1}
= (\hat f_2\circ \Psi_l^{-1})\circ (\Psi_l\circ u_\eps\circ \vphi^{-1})
\end{equation}
yielding the claim (as in the proof of Proposition \ref{emchar} (a)$\Rightarrow(b)$).

\noindent (c)$\Rightarrow$(b) follows from Proposition \ref{propneu}.

\noindent
(b)$\Rightarrow$(a): (i) in 
Definition \ref{hybridmod} being evident, let us
establish (ii). Let $L\comp (V,\vphi)$ a chart in $X$ and $L'\comp
(W,\Psi)$ a vector bundle chart in $F$ and choose an open
neighborhood $U$ of $L'$ whose compact closure is contained in $W$. Choose any
$l\in \{1,\dots, m\}$ and let $\hat f_l\in \mathrm{Hom}_c(F,\R\times\R^{m'})$
such that 
$$
\hat f_l|_{\pi_Y^{-1}(U)} = (\pro_l \times \mathrm{id}_{\R^{m'}}) 
\circ \Psi|_{\pi_Y^{-1}(U)}
$$
where $\pro_l:\R^m\to \R$. Let $p\in L\cap \underline{u_\eps}^{-1}(L')$. Then 
in some neighborhood of $p$ we have
\begin{equation}\label{dingsdanew}
(\hat f_l
\circ u_\eps\circ\vphi^{-1}) = (\psi_l\circ \underline{u_\eps}\circ \vphi^{-1},
{\bfs \psi}\circ u_\eps \circ \vphi^{-1})
\end{equation}
from which the claim follows.
\ep\ms

In the following result, $\sim_{h0}$ denotes the equivalence relation on $\esm^h[X,F]$
defined by imposing 
Definition \ref{hybridequiv} (i) and (ii) for $k=0$.
\bt \label{hybridequivchar}
 Let $(u_\eps)_\eps,\, (v_\eps)_\eps \in \esm^h[X,F]$. The following statements are equivalent:
\begin{itemize}
\item[(i)] $(u_\eps)_\eps \sim_{h} (v_\eps)_\eps$.
\item[(ii)] $(u_\eps)_\eps \sim_{h0} (v_\eps)_\eps$.
\item[(iii)] $(\hat f\circ u_\eps)\sim_{h} (\hat f\circ v_\eps)_\eps$ in 
$\esm^{h}[X,\R\times \R^{m'}]$ for all $\hat f\in \mathrm{Hom}_c(F,\R\times \R^{m'})$. 
\item[(iv)] $(\hat f\circ u_\eps)\sim_{h} (\hat f\circ v_\eps)_\eps$ in 
$\esm^{h}[X,\R\times \R^{m'}]$ for all $\hat f\in \mathrm{Hom}(F,\R\times \R^{m'})$. 
\end{itemize}
\et

\pr (i)$\Rightarrow$(ii) is clear.

\noindent(ii)$\Rightarrow$(iv): Let $\hat f\in \mathrm{Hom}(F,\R\times \R^{m'})$
and $L\comp V$, $(V,\vphi)$ a chart in $X$. Using the same notations as in the proof
of Proposition \ref{emhybridchar}, to show (iv) we have
to establish the $\ns$-estimates of order $0$ on $L$ for  
$(\underline{\hat f}\circ \underline{u_\eps}-
\underline{\hat f}\circ \underline{v_\eps})_\eps$ and 
$(\hat f_2 \circ u_\eps-\hat f_2 \circ v_\eps)_\eps$. 
Let $L'\comp Y$ and $\eps_0>0$ such that 
$\underline{u_\eps}(L)\cup \underline{v_\eps}(L) \subseteq L'$ for 
$\eps<\eps_0$ and cover $L'$ by open sets
$W_l'$, $\overline{W_l'}\comp W_l$ ($1\le l\le s$) with $(W_l,\Psi_l)$
vector bundle charts in $F$. 
Since $(\underline{u_\eps})_\eps \sim_0 (\underline{v_\eps})_\eps$,
for each small $\eps$ and each $p\in L$ there exists 
$l\in \{1,\dots,s\}$ such that $\underline{u_\eps}(p)$ and $\underline{v_\eps}(p)$ 
are contained in $W_l'$ simultaneously. Hence the first of the above estimates  
follows from Theorem  \ref{equivchar}. The second one 
follows from (\ref{locformfu2new}), 
Definition \ref{hybridequiv} (ii) (for $k=0$) and Lemma 2.5 of \cite{gfvm}: 
the constants in the ensuing estimates
can be chosen independently of $\eps$ since the domain of $\hat f_2 \circ \Psi_l^{-1}$
is an open neighborhood of ${\bfs \psi_l}(\overline{W_l'}) \times \R^{m'}$.

\noindent(iv)$\Rightarrow$(iii): Obvious.

\noindent(iii)$\Rightarrow$(i) Property (i) of 
Definition \ref{hybridequiv} is immediate from Theorem \ref{equivchar}. Property (ii) is established 
using an $\hat f_l$ as in the proof of
Proposition \ref{emhybridchar}, (b)$\Rightarrow$(a), employing representations as in (\ref{dingsdanew})
for both  $\hat f_l\circ u_\eps\circ \vphi^{-1}$ and $\hat f_l\circ v_\eps\circ \vphi^{-1}$.
\ep\ms

As in the previous cases we now use the above results to gain a point value description of
hybrid generalized functions.

\bt \label{pvcharhybrid} 
Let $u=[(u_\eps)_\eps] \in \gs^h[X,F]$, $\tilde p = [(p_\eps)_\eps] \in \widetilde X_c$. Then
\begin{itemize}
\item[(i)] $u(\tilde p) := [(u_\eps(p_\eps))_\eps]$ is a well-defined element of $F_c^{\sim_{vb}}$. 
\item[(ii)] If $u,\, v \in \gs^h[X,F]$ then $u=v$ if and only if $u(\tilde p) = v(\tilde p)$ 
for all $\tilde p\in \widetilde X_c$.
\end{itemize}
\et

\pr (i) Since $(p_\eps)_\eps$ is compactly supported, vb-moderateness of $(u_\eps(p_\eps))_\eps$ 
follows immediately from c-boundedness of $(\underline u_\eps)_\eps$ and 
Definition \ref{hybridmod}. Suppose that $(p_\eps)_\eps \sim (p_\eps')_\eps$. To show that 
$(u_\eps(p_\eps))_\eps \sim_{vb} (u_\eps(p_\eps'))_\eps$, by Theorem \ref{equivvbchar} and 
Corollary \ref{equivvbcharpoint} we may without loss of generality suppose that $F=\R\times \R^{m'}$.
Choose $K\comp X$ such that $p_\eps,\, p_\eps' \in K$ for $\eps$ small, cover $K$ by charts
$(V_1,\vphi_1),\dots,(V_l,\vphi_l)$ and write $K=\bigcup_{i=1}^l L_i$ with $L_i\comp V_i$. Then choosing the
trivial vector bundle chart on $\R\times \R^{m'}$ the claim follows from estimates of the form
$$
|u_\eps(p_\eps) - u_\eps(p_\eps')| \le \sup_{q\in L_i}\|D(u_\eps\circ \vphi_i^{-1})(\vphi_i(q))\|
\,|\vphi_i(p_\eps) - \vphi_i(p_\eps')|\,,
$$ 
and moderateness of $(u_\eps)_\eps$. Finally, let $(u_\eps)_\eps \sim_h (u_\eps')_\eps$. 
Then $(u_\eps(p_\eps))_\eps\sim_{vb} (u_\eps'(p_\eps))_\eps$ follows from 
Definition \ref{hybridequiv} (ii). 

\noindent (ii) Necessity has been established in (i). Conversely, suppose that $(u_\eps)_\eps
\not\sim_h (v_\eps)_\eps$, i.e., $(u_\eps)_\eps \not\sim_{h0} (v_\eps)_\eps$ 
(Theorem \ref{hybridequivchar}).
Then either $(\underline{u_\eps})_\eps \not\sim_0 (\underline{v_\eps})_\eps$ in which case by 
\cite{gfvm}, Prop.\ 2.14 we obtain a $\tilde p \in \widetilde X_c$ with $\underline u(\tilde p)
\not = \underline v(\tilde p)$, and, consequently, $u(\tilde p) \not= v(\tilde p)$. Hence the only
remaining possibility is that property (ii) for $k=0$ of 
Definition \ref{hybridequiv} is
violated. Hence there exist $m\in \N$, $L\comp V$ for a chart $(V,\vphi)$ in $X$, $L'\comp W$ 
for a vector bundle chart $(W,\Psi)$ in $F$ and sequences $\eps_j <1/j$, $p_j \in L \cap
\underline{u_{\eps_j}}^{-1}(L') \cap \underline{v_{\eps_j}}^{-1}(L')$ with
$$
|{\bfs \psi}\circ u_\eps(p_j) - {\bfs \psi}\circ v_\eps(p_j)| > j\eps_j^m\,.
$$
We set $p_\eps := p_j$ for $\eps_{j+1} <\eps \le \eps_j$. Then $\tilde p := [(p_\eps)_\eps] \in 
\widetilde X_c$ and $u(\tilde p) \not = v(\tilde p)$.
\ep\ms

Finally, compositions can be carried out unrestrictedly:
\bt Let $u=[(u_\eps)_\eps] \in \gs[X,Y]$, $v=[(v_\eps)_\eps] \in \gs^h[Y,G]$ and 
$w=[(w_\eps)_\eps]\in \mathrm{Hom}_{\gs}[G,H]$. Then $v\circ u:= 
[(v_\eps\circ u_\eps)_\eps]$ and $w\circ v:= [(w_\eps\circ v_\eps)_\eps]$
are well-defined elements of $\gs^h[X,G]$ and $\gs^h[Y,H]$, respectively.
\et

\pr Moderateness of $(v_\eps\circ u_\eps)_\eps$ and $(w_\eps\circ v_\eps)_\eps$ follows
in a straightforward way from the definitions. The respective classes are well-defined
(i.e., independent of the representatives of $u$, $v$, $w$) by Theorems \ref{pvchar},
\ref{pvcharhom} and \ref{pvcharhybrid} (ii).
\ep\ms

We next demonstrate that for any $u \in \gs[X,Y]$ the space 
$\mathrm{Hom}_u(E,$ $F):=\{v\in \mathrm{Hom}(E,$ $F) | \underline{v} = u\}$ can naturally be
endowed with the structure of an $\R$-vector space. Although the corresponding statement in the
smooth case is self-evident, in general the representatives 
$(v_\eps)_\eps$, $(v_\eps')_\eps$ of elements $v$, $v'$
of $\mathrm{Hom}_u(E,F)$ need not project onto the same representative $(u_\eps)_\eps$
of $u=\underline{v}=\underline{v'}\in\gs[X,Y]$, so that
simple fiberwise addition is in general not possible.
In order to handle this problem we employ the following result

\bp \label{fibertransplanting} 
Let $u=[(u_\eps)_\eps]\in \gs[X,Y]$ and $v\in \mathrm{Hom}_u(E,F)$. Then there
exists a representative $(v_\eps)_\eps$ of $v$ such that $\underline{v_\eps} = u_\eps$
for all $\eps\in I$.
\et

\pr We suppose that $X$ is connected for the moment, the modifications for the non-connected
case are then obvious. By $\mu_F: \R\times F \to F$ we denote fiberwise scalar multiplication.
Endow $Y$ with any Riemannian metric $h$ and cover $Y$ with $d_h$-balls $W_{\al}^{(r)}$ of radius 
$r$ such that $(W_\al^{(2r)}, \Psi_\al)$ is a vector bundle atlas of $F$. Choose a partition of
unity $\chi_j$ ($j\in \N$) subordinate to the covering $(W_{\al}^{(r)})_\al$ of $Y$ and let 
$\mathrm{supp} \chi_j \comp W_{\al_j}^{(r)}$ for $j\in \N$. Let $(\tilde v_\eps)_\eps$ be any
representative of $v$. For any relatively compact open subset $V$ of $X$ there exists some $\eps_0$
such that $\sup_{p\in V}d_h(\underline{\tilde v_\eps}(p),u_\eps(p)) < r$. 
Hence for $\eps<\eps_0$ and  $p\in V$, $u_\eps(p)\in\mathrm{supp}(\chi_j)$ implies
$\pi_Y(\tilde v_\eps(p))\in\mathrm{dom}(\psi_{\alpha_j})$ and we may define for
$e \in \pi_X^{-1}(V)$:
\begin{equation} \label{vlocdef}
v_\eps(e) := \sum_{j\in\N} \mu_F\left(\chi_j\circ u_\eps\circ \pi_X(e),\Psi_{\al_j}^{-1}
\left(\psi_{\al_j}(u_\eps(\pi_X(e))),{\bfs\psi}_{\al_j}(\tilde v_\eps(e))\right)\right)
\end{equation}
Since 
\beas
\tilde v_\eps(e) &=& \sum_{j\in \N} \mu_F(\chi_j\circ \underline{\tilde v_\eps}
\circ \pi_X(e),\tilde v_\eps(e)) \\
&=& \sum_{j\in\N} \mu_F\left(\chi_j\circ \underline{\tilde v_\eps}
\circ \pi_X(e),\Psi_{\al_j}^{-1}
\left(\psi_{\al_j}(\underline{\tilde v_\eps}(\pi_X(e))),
{\bfs\psi}_{\al_j}(\tilde v_\eps(e))\right)\right)
\eeas
it follows from $v\in \mathrm{Hom}_u(E,F)$ and the definition of $\sim_{vb}$ that
$(v_\eps)_\eps$ is a representative of $v|_V$. Finally, to construct a global 
representative of $v$ with the claimed properties we first note that by our
connectedness assumption $X$ is $\sigma$-compact, hence we may choose an exhaustive
sequence $K_m$ ($m\in\N$) of compact subsets of $X$ such that $K_m\comp
K_{m+1}^\circ$ for all $m$. By the above for each $m\in \N$ there exists $\eps_m>0$
such that $v_\eps$ as in (\ref{vlocdef}) is defined on $K_m^\circ$ for $\eps<\eps_m$. 
According to \cite{found}, Lemma 10.3 there exists a smooth $f: X\to (0,1]$ 
such that $0<f(p)\le \eps_n$ for $p\in K_n\setminus K_{n-1}^\circ$. Thus by 
defining $v$ according to (\ref{vlocdef}) for $(e,\eps) \in \{(e,\eps) \in E\times I |
\eps < f(\pi_X(e))\}$ (and arbitrarily elsewhere) we obtain the desired global 
representative of $v$.\ep

\bc For any $u\in \gs[X,Y]$, $\mathrm{Hom}_u(E,F)$ is a vector space.
\et

\pr Fixing any representative $(u_\eps)_\eps$ of $u$, for any $v_1$, $v_2 \in 
\mathrm{Hom}_u(E,F)$ we may choose representatives as in Proposition 
\ref{fibertransplanting}, so the claim follows.
\ep\ms

By a similar reasoning it follows that in case $X = Y$ the space of {\em strict} generalized 
vector bundle homomorphisms $\mathrm{Hom}_{id_X}(E,F)$ is naturally endowed 
with a vector space structure.

We finish this section demonstrating one of the consequences of the above result in the context
of the induced covariant derivative of generalized vector fields on a curve introduced in
\cite{curvature}, Sec.\ 5. Let $\al\in\gs[J,X]$, where $J$ is an interval in $\R$ and let 
$\hat g$ be a generalized metric on $X$ (cf.\ Def.\ 3.4 in \cite{curvature}) with generalized
Levi Civita connection $\hat D$ (cf.\ \cite{curvature}, Def.\ 5.1 and Th.\ 5.2). 
For any generalized vector field $\xi$ on $\al$, i.e., $\xi=[(\xi_\eps)_\eps]\in{\frak X}_{\gs}(\al)$
the induced covariant derivative $\xi'$ is defined componentwise, i.e., $\xi':=[(\xi'_\eps)_\eps]$,
where $\xi'_\eps$ is the classical induced covariant derivative with respect
to a representative $(\hat g_\eps)_\eps$ of $\hat g$. Note that by \cite{curvature}, Th.\ 3.1
any generalized pseudo-Riemannian metric (on any relatively compact open set) has a representative 
which consists of classical pseudo-Riemannian metrics for $\eps$ small enough. Now for $\xi$, $\eta\in
{\frak X}_{\gs}(\al)$ we may define $\hat g(\xi,\eta)$ componentwise since by Proposition
\ref{fibertransplanting} we may choose representatives $(\xi_\eps)_\eps$ and $(\eta_\eps)_\eps$ 
of $\xi$ respectively $\eta$ with $\underline{\xi_\eps}=\underline{\eta_\eps}=\al_\eps$ 
for some representative $(\al_\eps)_\eps$ of $\al$. Hence both sides of the equation
\beq\label{icd}
   \frac{d}{dt}\hat g(\xi,\eta)=\hat g(\xi',\eta)+\hat g(\xi,\eta')
\eeq
are well-defined. Moreover, equation (\ref{icd}) indeed holds by an application of the
classical result (cf.\ e.g. \cite{oneill}, Prop 3.18) at the level of
representatives.

\section{Coupled Calculus}\label{coupledcal}

In all variants of Colombeau generalized functions taking values in a linear space
compatibility with respect to the distributional setting is effected through 
the notion of {\em association}. We say an element $(u_\eps)_\eps$ of
$\esm(X)$ is associated to zero if $\int u_\eps\nu\to 0$ for all 
compactly supported one-densities $\nu$ on $X$ (see \cite{RD}, Sec.\ 2),
i.e., if $u_\eps \to 0$ weakly. 
If $(u_\eps)_\eps \in \ns(X)$ then it is automatically associated to
zero, so we may define a generalized function to be associated to zero
if its representatives are. We say that two generalized functions are
associated, $u\approx v$ if $u-v\approx 0$.
The equivalence relation $\approx$
allows us to define the linear quotient space
$\gs(X)/\approx$. Equality in $\gs(X)$ clearly implies association but
not conversely (e.g., any power $H^n$ ($n\ge2$) of the Heaviside function 
$H$ in $\gs(\R)$ is associated but not equal to $H$).

We say that $u \in \gs(X)$ admits $w \in \D'(X)$ as a distributional
shadow if
$$\lim_{\eps \to 0} \int u_\eps\nu=w(\nu), \quad 
\forall \nu \in \Gamma_c(\Vol(X)) $$

The distributional shadow of $u$ is uniquely determined (if
it exists).
Thus although the embedding $\iota$ of
$\D'(X)$ into $\gs(X)$ is not unique, one recovers a unique
description of $\D'(X)$ within the Colombeau algebra by working with 
$\iota(D'(X))/\approx$.

To extend the concept of association to the context of generalized functions
taking values in a smooth manifold we would need a ``diffeomorphism
invariant characterization'' of weak convergence. Based on the results
of this paper it is natural to aim at a notion of association defined by
$[(u_\eps)_\eps]\approx [(v_\eps)_\eps]\in\gs[X,Y]:\Leftrightarrow\  
[(f\circ u_\eps-f\circ v_\eps)_\eps]\approx 0\in\gs(X)\
\forall f\in\CC^\infty(Y)$.
However, such a construction does not reproduce the standard concept of 
association in the case that $Y$ is a vector space since
weak convergence cannot be
characterized by composition with smooth functions. To see this 
we first note that by \cite{MObook}, Ex.\ 10.6, we may choose $d_1,d_2\in \gs(\R)$,
both associated to $\delta$ such that $d_i^2\approx c_i\delta$ ($i=1,2$) with
$c_i\in\C$ and $c_1\not= c_2$. Hence $d_1\approx d_2$ but $f\circ d_1 \not\approx f\circ d_2$ for
$f= x \mapsto x^2$.


On the other hand the notion of $k$-association in $\gs(X)$ as introduced 
in \cite{ndg}, Sec.\ 5 can be generalized to the present setting using 
composition with smooth functions. 

\bd
\begin{itemize}
\item [(i)] Let $u=[(u_\eps)_\eps],\ v=[(v_\eps)_\eps]\in\gs[X,Y]$. 
We call $u$ and $v$ 
$k$-associated ($0\leq k\leq\infty$) with each other ($u\approx_k v$) if
$f\circ u\approx_k f\circ v\in \gs(X)$ for all $f\in\cinfty(Y)$.
\item [(ii)] We say that $u=[(u_\eps)_\eps]\in\gs[X,Y]$ admits 
$g:X\to Y$ as $\CC^k$-associated function ($u\approx_k g$) if 
$f\circ u\approx_k f\circ g\in\gs(X)$  for all $f\in\cinfty(Y)$.
\end{itemize}
\et

Analogously to the case of $\gs(X)$ \cite{ndg}, Def.\ 1 (ii) 
if $u \in \gs[X,Y]$ is $k$-associated to a function $g$, then 
$g\in\CC^k(X,Y)$. Furthermore if $u$ admits a $k$-associated
function at all, then the latter is unique.

If $Y=\R^n$ the notion defined in (i) above coincides with the one 
of \cite{ndg}, Def.\ 1(i).
To see this take
$[(u_\eps)_\eps]\approx_k[(v_\eps)_\eps]\in\gs(X)^n$. 
Since $\|u_\eps-v_\eps\|\to 0$ on compact sets the estimates on 
$\|f(u_\eps(x))-f(v_\eps(x))\|$
follow from \cite{gfvm}, Lemma 5. 
The derivatives can be estimated in the same way.
We have the following characterization of $0$-equivalence.
\bp \label{0asschar}
Let $(u_\eps)_\eps$, $(v_\eps)_\eps \in \gs[X,Y]$. The following statements
are equivalent:
\begin{itemize}
\item[(i)] $(u_\eps)_\eps \approx_0 (v_\eps)_\eps$.
\item[(ii)] For every Riemannian metric $h$ on $Y$ and every $K\comp X$,
$$
\sup_{p\in K} d_h(u_\eps(p),v_\eps(p)) \to 0 \qquad (\eps\to 0).
$$
\end{itemize}
\et
\pr This follows in complete analogy to the proof of the equivalence of (iii) and (iv)
of Theorem \ref{equivchar}.\ep



\bexs
\begin{itemize}
\item[(i)] By Theorem \ref{equivchar} two elements $(u_\eps)_\eps, (v_\eps)_\eps$ of 
$\gs[\R,\R]$ are equivalent if and only if for each $K\comp \R$ $\forall m\in \N$
$\sup_{x\in K} |u_\eps(x)-v_\eps(x)| = O(\eps^m)$. Hence the classes of $u_\eps(x) = \eps x$,
$v_\eps(x)=\eps^2x^2$ in $\gs[\R,\R]$ are different, but are easily seen to be
$0$-associated. 
\item[(ii)] In \cite{geo2}, solutions of the geodesic equation for a singular spacetime metric
with line-element
\begin{equation}
  \label{metric}
ds^2 = \delta(u)f(x,y) du^2 -du\, dv +dx^2 + dy^2  
\end{equation}
are calculated. Based on Theorem 1 of that paper it is shown in the final section of
\cite{curvature} that the $x$-component of the unique
solution of this ODE is c-bounded, hence can be regarded as a manifold-valued generalized function.
Furthermore, from \cite{geo2}, Theorem 3 it follows 
that it converges locally uniformly (hence by Proposition \ref{0asschar} is $0$-associated)
to a kink function. This corresponds to the physical expectation of broken geodesics of the
singular metric (\ref{metric}).
\end{itemize}
\et





\end{document}